\documentclass[11pt]{article}
\usepackage{amssymb}
\begin{document}
\setlength{\oddsidemargin}{0cm}
\setlength{\evensidemargin}{0cm}
\baselineskip=20pt

\begin{center} {\Large\bf
A Unified Algebraic Approach to Classical Yang-Baxter Equation }
\end{center}

\bigskip

\begin{center}  { \large Chengming ${\rm Bai}^{1,2}$} \end{center}

\begin{center}{\it 1. Chern Institute of
Mathematics \& LPMC, Nankai University, Tianjin 300071, P.R. China }
\end{center}

\begin{center}{\it 2. Dept. of Mathematics, Rutgers, The State University of
New Jersey, Piscataway, NJ 08854, U.S.A.}\end{center}

\vspace{0.3cm}

\begin{center} {\large\bf   Abstract } \end{center}

In this paper, the different operator forms of classical Yang-Baxter
equation are given in the tensor expression through a unified
algebraic method. It is closely related to left-symmetric algebras
which play an important role in many fields in mathematics and
mathematical physics. By studying the relations between
left-symmetric algebras and classical Yang-Baxter equation, we can
construct left-symmetric algebras from certain classical
$r$-matrices and conversely, there is a natural classical $r$-matrix
constructed from a left-symmetric algebra which corresponds to a
parak\"ahler structure in geometry. Moreover, the former in a
special case gives an algebraic interpretation of the
``left-symmetry'' as a Lie bracket ``left-twisted'' by a classical
$r$-matrix.

\vspace{0.2cm}



\bigskip

\section{Introduction}

Classical Yang-Baxter equation (CYBE) first arose in the study of
inverse scattering theory ([1-2]). It is also a special case of the
Schouten bracket in differential geometry which was introduced in
1940 ([3]). It can be regarded as a ``classical limit" of quantum
Yang-Baxter equation ([4]). They play a crucial role in many fields
like symplectic geometry, integrable systems, quantum groups,
quantum field theory and so on ([5] and the references therein).
Yang-Baxter system has become an important topic in both mathematics
and mathematical physics since 1980s.

The standard form of the CYBE in a Lie algebra is given in the
tensor expression as follows. Let ${\cal G}$ be a Lie algebra and
$r\in {\cal G}\otimes {\cal G}$. $r$ is called a solution of CYBE
in ${\cal G}$ if
$$[r_{12},r_{13}]+[r_{12},r_{23}]+[r_{13},r_{23}]=0\;\;{\rm in}\;U({\cal
G}), \eqno (1.1)$$ where $U({\cal G})$ is the universal enveloping
algebra of ${\cal G}$ and for $r=\sum\limits_i a_i\otimes b_i$,
$$r_{12}=\sum\limits_i a_i\otimes b_i\otimes 1;\;\; r_{13}=\sum\limits_i
a_i\otimes1\otimes b_i ;\;\;r_{23}=\sum\limits_i 1\otimes a_i\otimes
b_i.\eqno (1.2)$$ $r$ is also called a classical $r$-matrix due to
the expression of $r$ under a basis of ${\cal G}$.

There are a lot of results on CYBE when ${\cal G}$ is semisimple
(cf. [6], [7], etc.). However, it is not easy to study equation
(1.1) directly in a general case. A natural idea is to replace the
tensor form by a linear operator. There are several approaches. In
[8], Semonov-Tian-Shansky studied CYBE systematically. In
particular, an operator form of CYBE is given as a linear map
$R:{\cal G}\rightarrow {\cal G}$ satisfying
$$[R(x), R(y)]=R([R(x),y]+[x,R(y)]),\;\;\forall x,y\in {\cal
G}.\eqno (1.3)$$ It is equivalent to the tensor form  (1.1) of CYBE
when the following two conditions are satisfied: (a) there exists a
nondegenerate symmetric invariant bilinear form on ${\cal G}$ and
(b) $r$ is skew-symmetric. However, the relation between the
operator form (1.3) and the tensor form (1.1) in a general case is
still not clear.

Later, Kupershmidt re-studied the CYBE in [9]. When $r$ is
skew-symmetric, the tensor form (1.1) of CYBE is equivalent to a
linear map $r: {\cal G}^*\rightarrow {\cal G}$ satisfying
$$[r(x),r(y)]=r({\rm ad}^* r(x)(y)-{\rm
ad}^*r(y)(x)),\;\forall x,y\in {\cal G}^*,\eqno (1.4)$$ where ${\cal
G}^*$ is the dual space of ${\cal G}$ and ${\rm ad}^*$ is the dual
representation of adjoint representation (coadjoint representation)
of the Lie algebra ${\cal G}$. Moreover, Kupershmidt generalized the
above ${\rm ad}^*$ to be an arbitrary representation $\rho:{\cal
G}\rightarrow gl(V)$ of ${\cal G}$, that is, a linear map
$T:V\rightarrow {\cal G}$ satisfying
$$[T(u), T(v)]=T(\rho(T(u))v-\rho(T(v))u),\forall u,v\in V,\eqno (1.5)$$
which was regarded as a natural generalization of CYBE. Such an
operator is called an ${\cal O}$-operator associated to $\rho$. Note
that the operator form (1.3) of CYBE given by Semonov-Tian-Shansky
is just an ${\cal O}$-operator associated to the adjoint
representation of ${\cal G}$. However, there is not a direct
relation between the ${\cal O}$-operators and the tensor form (1.1)
of CYBE, either.

In this paper, we give a further study on CYBE which unifies the
above different operator forms of CYBE. The key is how to interpret
the ${\cal O}$-operators in terms of the tensor expression. Our idea
is to extend the Lie algebra ${\cal G}$ to construct a bigger Lie
algebra such that the ${\cal O}$-operators can be related to the
solutions of the tensor form of CYBE in it. Thus we can obtain not
only the direct relations between the above operators and the tensor
form (1.1) of CYBE, but also the corresponding results of
Semonov-Tian-Shansky and Kupershmidt as special cases. It is quite
similar to the double construction ([7]).

Furthermore, there is an algebraic structure behind the above study.
It is the left-symmetric algebra (or under other names like pre-Lie
algebra, quasi-associative algebra, Vinberg algebra and so on).
Left-symmetric algebras are a class of nonassociative algebras
coming from the study of convex homogeneous cones, affine manifolds
and affine structures on Lie groups, deformation of associative
algebras ([10-13]) and then appear in many fields in mathematics and
mathematical physics, such as complex and symplectic structures on
Lie groups and Lie algebras ([14-18]), integrable systems ([19-20]),
Poisson brackets and infinite-dimensional Lie algebras ([21-23]),
vertex algebras ([24]), quantum field theory ([25]), operads ([26])
and so on (more examples can be found in a survey in [27] and the
references therein).

Although some scattered results are known in certain references (cf.
[9,28-31], etc.), we give a systematic study on the relations
between left-symmetric algebras and CYBE in this paper. It can be
regarded as a generalization of the correspondence between
left-symmetric algebras and bijective 1-cocycles whose existence
gives a necessary and sufficient condition for a Lie algebra with a
compatible left-symmetric algebra structure. For the study of
left-symmetric algebras, it provides a construction from certain
classical $r$-matrices. In particular, as a special case, an
algebraic interpretation of the so-called ``left-symmetry'' is
induced by equation (1.3): in certain sense, the ``left-symmetry"
can be interpreted as a Lie bracket ``left-twisted'' by a classical
$r$-matrix. On the other hand, for the study of CYBE, there is a
natural classical $r$-matrix with a simple form constructed from a
left-symmetric algebra which corresponds to a parak\"ahler structure
in geometry.

The paper is organized as follows. In Section 2, we construct a
direct relation between ${\cal O}$-operators and the tensor form of
CYBE. In the cases of adjoint representations and co-adjoint
representations, we can get the operator forms (1.3) and (1.4) of
CYBE. In section 3, we briefly introduce left-symmetric algebras and
then study the relations between them and CYBE. In section 4, we
summarize the main results obtained in the previous sections.

Throughout this paper, without special saying, all algebras are of
finite dimension and over an algebraically closed field of
characteristic 0 and $r$ is a solution of CYBE or $r$ is a classical
$r$-matrix refers to that $r$ satisfies the tensor form (1.1) of
CYBE.

\section{ ${\cal O}$-operators and the tensor form of CYBE}

At first, we give some notations. Let ${\cal G}$ be a Lie algebra
 and $r\in {\cal G}\otimes {\cal G}$. $r$ is said to be skew-symmetric if
$$r=\sum\limits_i (a_i\otimes b_i-b_i\otimes a_i).\eqno (2.1)$$
For $r=\sum\limits_i a_i\otimes b_i \in {\cal G}\otimes {\cal G}$,
we denote
$$r^{21}=\sum\limits_i b_i\otimes a_i.\eqno (2.2)$$
On the other hand, let $\rho:{\cal G}\rightarrow gl(V)$ be a
representation of the Lie algebra ${\cal G}$. On the vector space
${\cal G}\oplus V$, there is a natural Lie algebra structure
(denoted by ${\cal G}\ltimes_\rho V$) given as follows ([32]).
$$[x_1+v_1, x_2+v_2]=[x_1,x_2]+\rho (x_1)v_2-\rho (x_2)v_1,\;\forall
x_1,x_2\in {\cal G}, v_1, v_2\in V. \eqno (2.3)$$

Let $\rho^*:{\cal G}\rightarrow gl(V^*)$ be the dual representation
of the representation $\rho: {\cal G}\rightarrow gl(V)$ of the Lie
algebra ${\cal G}$. Then there is a close relation between the
${\cal O}$-operator associated to $\rho$ and the (skew-symmetric)
solutions of CYBE in ${\cal G}\ltimes_{\rho^*}V^*$.

Any linear map $T:V\rightarrow {\cal G}$ can be identified as an
element in ${\cal G}\otimes V^*\subset ({\cal
G}\ltimes_{\rho^*}V^*)\otimes ({\cal G}\ltimes_{\rho^*}V^*)$ as
follows. Let $\{e_1,\cdots,e_n\}$ be a basis of ${\cal G}$. Let
$\{v_1,\cdots, v_m\}$ be a basis of $V$ and $\{ v_1^*,\cdots,
v_m^*\}$ be its dual basis, that is $v_i^*(v_j)=\delta_{ij}$. Set
$T(v_i)=\sum\limits_{j=1}^na_{ij}e_j, i=1,\cdots, m$. Since as
vector spaces, ${\rm Hom}(V,{\cal G})\cong {\cal G}\otimes V^*$, we
have
$$T=\sum_{i=1}^m T(v_i)\otimes v_i^*=\sum_{i=1}^m\sum_{j=1}^n
a_{ij}e_j\otimes v_i^*\in {\cal G}\otimes V^*\subset ({\cal
G}\ltimes_{\rho^*}V^*)\otimes ({\cal G}\ltimes_{\rho^*}V^*).\eqno
(2.4)$$

\noindent {\bf Claim}\quad {\it $r=T-T^{21}$ is a skew-symmetric
solution of CYBE in ${\cal G}\ltimes_{\rho^*}V^*$ if and only if $T$
is an ${\cal O}$-operator.

In fact, by equation (2.4), we have

{\small\begin{eqnarray*} &&[r_{12},r_{13}] =\sum_{i,k=1}^m
\{[T(v_i), T(v_k)]\otimes v_i^*\otimes v_k^* -\rho^*
(T(v_i))v_k^*\otimes v_i^*\otimes T(v_k)+ \rho^*
(T(v_k))v_i^*\otimes
T(v_i)\otimes v_k^* \};\\
&&[ r_{12},r_{23}] =\sum_{i,k=1}^m \{- v_i^*\otimes [T(v_i), T(v_k)]\otimes
v_k^*
- T(v_i)\otimes \rho^* (T(v_k))v_i^*\otimes v_k^*+ u_i^*\otimes \rho^*
(T(v_i))v_k^*\otimes
T(v_k)\};\\
&&[r_{13},r_{23}] =\sum_{i,k=1}^m \{v_i^*\otimes u_k^*\otimes [T(v_i),
T(v_k)]
+ T(v_i)\otimes v_k^*\otimes \rho^*(T(v_k))v_i^*- v_i^*\otimes T(v_k)\otimes
\rho^* (T(v_i))v_k^*\}.
\end{eqnarray*}}
By the definition of dual representation, we know
$$\rho^*(T(v_k))v_i^*=-\sum_{j=1}^m v_i^*(\rho(T(v_k))v_j) v_j^*.$$
Thus
\begin{eqnarray*}
&&-\sum_{i,k=1}^m T(v_i)\otimes \rho^*(T(v_k))v_i^*\otimes
v_k^*=-\sum_{i,k=1}^mT(v_i)\otimes
[\sum_{j=1}^m -v_i^*(\rho(T(v_k))v_j) v_j^*]\otimes v_k^*\\
&&=\sum_{i,k=1}^m \sum_{j=1}^m v_j^*(\rho(T(v_k))v_i) T(v_j)\otimes
v_i^*\otimes v_k^*
=\sum_{i,k=1}^m T(\sum_{j=1}^m (v_j^*(\rho(T(v_k))v_i) v_j)\otimes
v_i^*\otimes v_k^* \\
&&=\sum_{i,k=1}^m T(\rho (T(v_k))v_i)\otimes v_i^*\otimes v_k^*.
\end{eqnarray*}
Therefore
\begin{eqnarray*}
&&[r_{12},r_{13}]+[r_{12},r_{23}]+[r_{13},r_{23}]\\
&&=\sum_{i,k=1}^m\{
([T(v_i,v_k)]+T(\rho(T(v_k))v_i)-T(\rho(T(v_i))v_k))\otimes v_i^*\otimes
v_k^*\\
&& -v_i^*\otimes ([T(v_i,v_k)]+T(\rho(T(v_k))v_i)-T(\rho(T(v_i))v_k))\otimes
v_k^*\\
&& + v_i^*\otimes v_k^*\otimes
([T(v_i,v_k)]+T(\rho(T(v_k))v_i)-T(\rho(T(v_i))v_k))\}.
\end{eqnarray*}
So $r$ is a classical $r$-matrix in  ${\cal G}\ltimes_{\rho^*}V^*$
if and only if $T$ is an ${\cal O}$-operator.}

Obviously, the above $r$ is exactly the skew-symmetric classical
$r$-matrix in ${\cal G}\ltimes_{\rho^*}V^*$ which is in ${\cal
G}\otimes V^*- V^*\otimes {\cal G}$.

Next we consider the cases that $\rho$ is the adjoint representation
${\rm ad}:{\cal G}\rightarrow gl({\cal G})$ or the coadjoint
representation ${\rm ad}^*:{\cal G} \rightarrow gl({\cal G}^*)$ with
$<{\rm ad}^*x(y^*),z>=-<y^*,[x,z]>$ for any $x,z\in {\cal G}$ and
$y^*\in {\cal G}^*$, where $<,>$ is the ordinary pair between ${\cal
G}$ and ${\cal G}^*$.

{\bf Case (I)}\quad $\rho={\rm ad}^*$, the coadjoint representation.
In this case $V={\cal G}^*$ and $V^*={\cal G}$. For any linear map
$T:{\cal G}^*\rightarrow {\cal G}$, $T$ can be identified as an
element in ${\cal G}\otimes {\cal G}$ by
$$<T(u),v>=<u\otimes v, T>,\;\;\forall\; u,v\in {\cal G}^*.\eqno
(2.5)$$ Therefore, although $r=T-T^{21}\in ({\cal G}\ltimes_{\rm
ad}{\cal G})\otimes ({\cal G}\ltimes_{\rm ad}{\cal G})$, in fact,
$r\in {\cal G}\otimes {\cal G}$. Hence $r$ is a skew-symmetric
solution of CYBE in ${\cal G}$ if and only if $T$ satisfies equation
(1.4).

In particular, suppose that $r\in {\cal G}\otimes {\cal G}$ which is
identified as a linear map from ${\cal G}^*$ to ${\cal G}$ and $r$
itself is skew-symmetric. Then $r-r^{21}=2r$. Obviously $r$ is a
classical $r$-matrix if and only if $2r$ is also a classical
$r$-matrix. Therefore $r$ is a solution of CYBE in ${\cal G}$ if and
only if $r$ satisfies equation (1.4). Thus, in this case, the CYBE
in the tensor expression (1.1) is equivalent to equation (1.4),
which was given by Kupershmidt ([9]).

{\bf Case (II)}\quad $\rho={\rm ad}$, the adjoint representation. In
this case, we suppose that the Lie algebra ${\cal G}$ is equipped
with a nondegenerate invariant symmetric bilinear form  $B(\;,\;)$.
That is,
$$B(x,y)=B(y,x),\;\;B([x,y],z)=B(x,[y,z]),\;\forall x,y,z\in {\cal G}.
\eqno (2.6)$$ Hence ${\cal G}^*$ is identified with ${\cal G}$. Let
$r\in {\cal G}\otimes {\cal G}$. Then $r$ can be identified as a
linear map from ${\cal G}$ to ${\cal G}$. If $r$ is skew-symmetric,
then $r$ is a solution of CYBE  in ${\cal G}$ if and only if $r$
satisfies equation (1.3). Thus, in this case, the CYBE in the tensor
expression (1.1) is equivalent to equation (1.3), which was given by
Semonov-Tian-Shansky ([8]).

\section{CYBE and left-symmetric algebras}

A left-symmetric algebra $A$ is a vector space over a field ${\bf
F}$ equipped with a bilinear product $(x,y)\rightarrow xy$
satisfying that for any $x,y,z\in A$, the associator
$$(x,y,z)=(xy)z-x(yz)\eqno (3.1)$$
is symmetric in $x,y$, that is,
$$(x,y,z)=(y,x,z),\;\;{\rm or}\;\;{\rm
equivalently}\;\;(xy)z-x(yz)=(yx)z-y(xz).\eqno (3.2)$$

Left-symmetric algebras are Lie-admissible algebras (cf. [33-34]).
In fact, let $A$ be a left-symmetric algebra.  Then the commutator
$$[x,y]=xy-yx,\;\;\forall x,y\in A,\eqno (3.3)$$
defines a Lie algebra ${\cal G}(A)$, which is called the
sub-adjacent Lie algebra of $A$ and $A$ is also called the
compatible left-symmetric algebra structure on the Lie algebra
${\cal G}(A)$.

Furthermore, for any $x\in A$, let $L_x$ denote the left
multiplication operator, that is, $L_x(y)=xy$ for any $y\in A$. Then
$L:{\cal G}(A)\rightarrow gl({\cal G}(A))$ with $x\rightarrow L_x$
gives a regular representation of the Lie algebra ${\cal G}(A)$,
that is,
$$[L_x,L_y]=L_{[x,y]},\;\;\forall x,y\in A. \eqno (3.4)$$

It is not true that there is a compatible left-symmetric algebra
structure on every Lie algebra. For example, a real or complex Lie
algebra ${\cal G}$ with a compatible left-symmetric algebra
structure must satisfy the condition $[{\cal G}, {\cal G}]\ne {\cal
G}$ ([34]), hence there does not exist a compatible left-symmetric
algebra structure on any real or complex semisimple Lie algebra.
Here, we briefly introduce a necessary and sufficient condition for
a Lie algebra with a compatible left-symmetric algebra structure
([33]). Let ${\cal G}$ be a Lie algebra and $\rho:{\cal
G}\rightarrow gl(V)$ be a representation of ${\cal G}$. A 1-cocycle
$q$ associated to $\rho$ (denoted by $(\rho,q)$) is defined as a
linear map from ${\cal G}$ to $V$ satisfying
$$q[x,y]=\rho(x)q(y)-\rho(y)q(x),\forall x,y\in {\cal G}.\eqno (3.5)$$
Then there is a compatible left-symmetric algebra structure on
${\cal G}$ if and only there exists a bijective 1-cocycle of ${\cal
G}$. In fact, let $(\rho,q)$ be a bijective 1-cocycle of ${\cal G}$,
then
$$x*y=q^{-1}\rho(x)q(y),\;\;\forall x,y\in {\cal G},\eqno (3.6)$$
defines a compatible left-symmetric algebra structure on ${\cal G}$.
Conversely, for a left-symmetric algebra $A$, $(L,id)$ is a
bijective 1-cocycle of ${\cal G}(A)$, where $id$ is the identity
transformation on ${\cal G}(A)$.

Note that for any Lie algebra ${\cal G}$ and its representation
$\rho:{\cal G}\rightarrow gl(V)$,  a linear isomorphism
$T:V\rightarrow {\cal G}$(hence $\dim {\cal G}=\dim V$) is an ${\cal
O}$-operator associated to $\rho$ if and only if $T^{-1}$ is a
(bijective) 1-cocycle of ${\cal G}$ associated to $\rho$. Therefore,
if $q:{\cal G}\rightarrow V$ is a bijective 1-cocycle of ${\cal G}$
associated to $\rho$, then $q^{-1}-(q^{-1})^{21}$ is a solution of
CYBE in ${\cal G}\ltimes_{\rho^*}V^*$. In particular, let $({\cal
G},*)$ be a left-symmetric algebra. Since $T=id$ is an ${\cal
O}$-operator associated to the regular representation $L$, we have
$$r=\sum_{i}^n (e_i\otimes e_i^*-e_i^*\otimes e_i)\eqno (3.7)$$
is a solution of CYBE in ${\cal G} \ltimes_{L^*} {\cal G}^*$, where
$\{e_1,\cdots, e_n\}$ is a basis of ${\cal G}$ and $\{e_1^*,\cdots,
e_n^*\}$ is its dual basis. Moreover, we would like to point out
that in the Lie algebra ${\cal G}(A)\ltimes_{{\rm ad}^*}{\cal
G}^*(A)$, equation (3.7) is not a solution of CYBE, but a solution
of the so-called modified CYBE ([8]). That is, in ${\cal
G}(A)\ltimes_{{\rm ad}^*}{\cal G}^*(A)$, equation (3.7) does not
satisfy equation (1.1), but satisfies
$$[x\otimes1\otimes 1+1\otimes x\otimes 1+1\otimes1\otimes x,\;\;
[r_{12},r_{13}]+[r_{12},r_{23}]+[r_{13},r_{23}]]=0,\forall x\in
{\cal G}(A)\ltimes_{{\rm ad}^*}{\cal G}^*(A).\eqno (3.8)$$

On the other hand, if an ${\cal O}$-operator $T$ associated to
$\rho$ is invertible, then $T^{-1}$ is a bijective 1-cocycle of
${\cal G}$ associated to $\rho$. Hence
$$x\cdot y=T(\rho (x) (T^{-1}(y))),\;\;\forall\; x,y\in {\cal G}. \eqno (3.9)$$
 defines a compatible left-symmetric algebra structure on ${\cal G}$ through
equation (3.6) by letting $q=T^{-1}$. Moreover, for any $u,v\in V$,
let $x=T(u), y=T(v)$. Thus by equation (3.9), we have
$$T(u)\cdot T(v)=T(\rho(T(u))v).$$
Since $T$ is invertible, there exists a left-symmetric algebra
structure on $V$ induced from the left-symmetric algebra structure
on ${\cal G}$ by
$$u*v=T^{-1}(T(u)\cdot T(v))=\rho(T(u))(v),\;\;\forall\; u,v\in V.\eqno (3.10)$$
It is obvious that $T$ is an isomorphism of left-symmetric algebras
between them

Furthermore, we can generalize the above construction of
left-symmetric algebras to a general ${\cal O}$-operator. Let ${\cal
G}$ be a Lie algebra and $\rho:{\cal G}\rightarrow gl(V)$ be its
representation.  Let $T:V\rightarrow {\cal G}$ be a linear map. Then
on $V$, the new product
$$u*v=\rho(T(u))v,\;\;\forall u,v\in V\eqno (3.11)$$
satisfies that for any $u,v,w\in V$,
\begin{eqnarray*}
(u,v,w)-(v,u,w)&=&\rho (T\rho (T(u))v)w-\rho(T(u))\rho(T(v))w
-(T\rho (T(v))u)w+\rho(T(v))\rho(T(u))w\\
&=&\rho ([T(v),T(u)])w+\rho (T(\rho(T(u))v-\rho(T(v))u))w.
\end{eqnarray*}
Hence equation (3.11) defines a left-symmetric algebra if and only
if
$$[T(u),T(v)]-T(\rho(T(u))v-\rho(T(v))u)\in {\rm Ker}\rho,\;\;\forall
u,v\in V,\eqno (3.12)$$ where ${\rm Ker}\rho=\{ x\in {\cal
G}|\rho(x)=0\}$. In particular, for any ${\cal O}$-operator
$T:V\rightarrow {\cal G}$ associated to $\rho$, equation (3.11)
defines a left-symmetric algebra on $V$. Therefore $V$ is a Lie
algebra as the sub-adjacent Lie algebra of this left-symmetric
algebra and $T$ is a Lie algebraic homomorphism. Furthermore,
$T(V)=\{T(v)|v\in V\}\subset {\cal G}$ is a Lie subalgebra of ${\cal
G}$ and there is an induced left-symmetric algebra structure on
$T(V)$ given by
$$T(u)\cdot T(v)=T(u*v),\;\;\forall u,v\in V.\eqno (3.13)$$
Moreover, its sub-adjacent Lie algebra structure is just the Lie
subalgebra structure of ${\cal G}$ and $T$ is a homomorphism of
left-symmetric algebras.

In fact, the above Lie algebra structure on $V$ coincides the
standard construction of Lie bialgebra from the classical $r$-matrix
$r=T-T^{21}$ as follows. ${\cal G}\ltimes_{\rho^*}V^*$ is a Lie
bialgebra ([5,7]) with the cobracket $\delta (f)=[f\otimes
1+1\otimes f, r]$ for any $f\in {\cal G}\ltimes_{\rho^*}V^*$. Hence
$V^*$ is a Lie co-subalgebra of ${\cal G}\ltimes_{\rho^*}V^*$.
Therefore $V$ is a Lie algebra just given by $[u,v]=\rho
(T(u))v-\rho(T(v))u$ for any $u,v\in V$.

According to Drinfel'd ([7]), $r\in {\cal G}\otimes {\cal G}$ is a
skew-symmetric and nondegenerate solution of CYBE in ${\cal G}$ if
and only if the bilinear form $B$ on ${\cal G}$ given by
$$B(x,y)=<r^{-1}(x),y>,\;\;\forall x,y\in {\cal G}, \eqno (3.14)$$
is a 2-cocycle on ${\cal G}$, that is,
$$B([x,y],z)+B([y,z],x)+B([z,x],y)=0,\;\;\forall x,y,z\in {\cal
G}.\eqno (3.15)$$ In geometry, a skew-symmetric and nondegenerate
2-cocycle on a Lie algebra ${\cal G}$ is also called a symplectic
form which corresponds to a symplectic form on a Lie group whose Lie
algebra is ${\cal G}$ and such a Lie group (or Lie algebra) is
called a symplectic Lie group (or Lie algebra) ([15]).

Now, let us return to our study on left-symmetric algebras and CYBE.
Let $A$ be a left-symmetric algebra. Then the classical $r$-matrix
$r$ given by equation (3.7) is nondegenerate. Moreover $r:({\cal G}
\ltimes_{L^*} {\cal G}^*)^*$ $\rightarrow$ ${\cal G} \ltimes_{L^*}
{\cal G}^*$ satisfies the following equations:
$$r(e_i^*)=e_i^*,\;\;r(e_i)=-e_i,\;\;i=1,\cdots,n.$$
Therefore, for any $i,j,k,l$, we have
$$<r^{-1}(e_i+e_j^*),e_k+e_l^*>=<-e_i+e_j^*,e_k+e_l^*>=<-e_i,e_l^*>+<e_k,e_j^*>.$$
Hence there is a natural 2-cocycle $\omega$ (symplectic form) on
${\cal G}\ltimes_{L^*} {\cal G}^*$ induced by $r^{-1}: {\cal G}
\ltimes_{L^*} {\cal G}^*\rightarrow ({\cal G} \ltimes_{L^*} {\cal
G}^*)^*$ given by
$$\omega(x+x^*,y+y^*)=<x^*,y>-<y^*,x>,\;\; \forall\; x,y\in {\cal G}, x^*,y^*\in {\cal
G}^*.\eqno (3.16)$$

The above structure ${\cal G} \ltimes_{L^*} {\cal G}^*$ with the
symplectic form $\omega$ (3.16) corresponds to a parak\"ahler
structure. In geometry, a parak\"ahler manifold is a symplectic
manifold with a pair of transversal Lagrangian foliations ([35]). A
parak\"ahler Lie algebra ${\cal G}$ is just the Lie algebra of a Lie
group $G$ with a $G$-invariant parak\"ahler structure ([36]). On the
other hand, such a structure is just a phase space of ${\cal G}$ in
mathematical physics ([37-39]).

Next we still consider the case that $\rho={\rm ad}$, the adjoint
representation. Let ${\cal G}$ be a Lie algebra and $f$ be a linear
transformation on ${\cal G}$. Then on ${\cal G}$ the new product
$$x*y=[f(x),y],\forall x,y\in {\cal G}\eqno (3.17)$$
defines a left-symmetric algebra if and only if
$$[f(x),f(y)]-f([f(x),y]+[x,f(y)])\in C({\cal G}),\;\forall x,y\in {\cal
G},\eqno (3.18)$$ where $C({\cal G})=\{ x\in {\cal
G}|[x,y]=0,\;\forall\; y\in {\cal G}\}$ is the center of Lie algebra
${\cal G}$. In particular, the map given by equation (1.3) defines a
left-symmetric algebra on ${\cal G}$ through equation (3.17). On the
other hand, if in addition, the center $C({\cal G})$ is zero, then
the linear transformation $f$ satisfying equation (3.18) just
satisfies equation (1.3), that is, $f$  satisfies the operator form
of CYBE.

The formula (3.17) was also given in [31] and a similar construction
for Novikov algebras (left-symmetric algebras with commutative right
multiplication operators) was given with $r$ satisfying some
additional conditions in [40]. We would like to point out that the
above construction cannot get all left-symmetric algebras.

Moreover, the above discussion gives an algebraic interpretation of
``left-symmetry" (3.2). Let $\{e_i\}$ be a basis of Lie algebra
${\cal G}$ and $r$ be a linear transformation satisfying equation
(1.3). Set $r(e_i)=\sum_{j\in I} r_{ij}e_j$. Then the
basis-interpretation of equation (3.17) is given as
$$e_i*e_j=\sum_{l\in I}r_{il}[e_l, e_j].\eqno (3.19)$$
In this sense, such a construction of left-symmetric algebras can be
regarded as a Lie algebra ``left-twisted'' by a classical
$r$-matrix. On the other hand, we consider the right-symmetric
algebra, that is, $(x,y,z)=(x,z,y)$ for any $x,y,z\in A$, where
$(x,y,z)$ is the associator given by equation (3.1). Set
$$e_i\cdot e_j=[e_i,r(e_j)]=\sum_{l\in I}r_{jl}[e_i, e_l].\eqno (3.20)$$
Then the above product defines a right-symmetric algebra on ${\cal
G}$, which can be regarded as a Lie algebra ``right-twisted'' by a
classical $r$-matrix.

Roughly speaking, the CYBE describes certain ``permutation''
relation and left-symmetry or right-symmetry is a kind of special
``permutation''. There is a close relation between them by equation
(3.19) or (3.20).

At the end of this section, we give a further study on the linear
transformations satisfying equation (3.18). For a left-symmetric
algebra structure on ${\cal G}$ given by equation (3.17) with $f$
satisfying equation (3.18), if in addition, its sub-adjacent Lie
algebra is just ${\cal G}$ itself, that is,
$$[x,y]=[f(x),y]-[y,f(x)],\;\;\forall x,y\in {\cal G},\eqno (3.21)$$
then it is a left-symmetric inner derivation algebra, that is, for
every $x\in {\cal G}$, $L_x$ is an interior derivation of the Lie
algebra ${\cal G}$. Such a structure corresponds to a flat
left-invariant connection adapted to the interior automorphism
structure of a Lie group, which was first studied in [34]. Moreover,
every left-symmetric inner derivation algebra can be obtained by
this way. On the other hand, if $f$ satisfies equation (1.3) and $f$
is invertible, then $f$ satisfies equation (3.21) if and only if $f$
is an automorphism of ${\cal G}$. Under this condition, ${\cal G}$
must be solvable since the sub-adjacent Lie algebra of a
left-symmetric inner derivation algebra is solvable ([34]).

Moreover, there is a general conclusion. Let ${\cal G}$ be a complex
Lie algebra with a nondegenerate symmetric invariant bilinear form.
If there exists an invertible skew-symmetric classical $r$-matrix
$r$, then $r$ can be identified as a linear transformation on ${\cal
G}$ satisfying equation (1.3). Hence $r^{-1}$ is a derivation of
${\cal G}$. Since a complex Lie algebra with a nondegenerate
derivation must be nilpotent (cf. [41]), ${\cal G}$ is nilpotent.
This conclusion also generalizes a similar result for complex
semisimple Lie algebras ([5], Proposition 2.2.5).

\section{Summary}

In this paper, we interpret  the ${\cal O}$-operators in terms of
the tensor expression. Thus the different operator forms of CYBE are
given in the tensor expression through a unified algebraic method.
Since the Lie bialgebra structures are obtained through the
solutions of CYBE in the tensor form as
$$\delta(x)=[x\otimes 1+1\otimes x, r],\;\;\forall x\in {\cal G},$$
it is easy to get the corresponding Lie bialgebra structures from
the different operator forms of CYBE through our study. It will be
also useful to consider the quantization of these Lie bialgebra
structures (for example, try to find the corresponding Drinfel'd
quantum twist [5,7]).

On the other hand, there are close relations between left-symmetric
algebras and CYBE. They can be regarded as a generalization of the
study of ${\cal O}$-operators in the cases that the ${\cal
O}$-operators are invertible. We can construct left-symmetric
algebras from certain classical $r$-matrices and conversely, there
is a natural classical $r$-matrix constructed from a left-symmetric
algebra which corresponds to a parak\"ahler structure in geometry.
Moreover, the former in a special case gives an algebraic
interpretation of the ``left-symmetry'' as a Lie bracket
``left-twisted'' by a classical $r$-matrix.

\begin{center}{\bf Acknowledgements}\end{center}

The author thanks Professors I.M. Gel'fand, B.A. Kupershmidt and P.
Etingof for useful suggestion and great encouragement. In
particular, the author is grateful for Professor P. Etingof who
invited the author to visit MIT and the work was begun there. The
author also thanks Professors J. Lepowsky, Y.Z. Huang and H.S. Li
for the hospitality extended to him during his stay at Rutgers, The
State University of New Jersey and for valuable discussions. This
work was supported in part by S.S. Chern Foundation for Mathematical
Research, the National Natural Science Foundation of China
(10571091, 10621101), NKBRPC (2006CB805905), Program for New Century
Excellent Talents in University and K.C. Wong Education Foundation.

\bigskip

\begin{center}{\bf Reference}\end{center}
\baselineskip=16pt
\begin{description}
\item[[1]] L.D. Faddeev, L. Takhtajan, The quantum inverse
scattering method of the inverse problem and the Heisenberg XYZ
model, Russ. Math. Surv. 34 (1979) 11-68.
\item[[2]] L.D. Faddeev, L. Takhtajan, Hamiltonian methods in the
theory of solitons, Springer, Berlin (1987).
\item[[3]] J.A. Schouten, \"Uber differentialkomitanten zweier
kontravarianter Gr\"ossen, Nede. Akad. Wete. Proc. A 43 (1940)
449-452.
\item[[4]] A.A. Belavin, Dynamical symmetry of integrable quantum
systems, Nucl. Phys. B 180 (1981) 189-200.
\item[[5]] V. Chari, A. Pressley, A guide to quantum groups, Cambridge
University Press, Cambridge (1994).
\item[[6]] A.A. Belavin, V.G. Drinfel'd, Solutions of classical Yang-Baxter
equation for simple Lie algebras, Funct. Anal. Appl. 16 (1982)
159-180.
\item[[7]] V. Drinfel'd, Hamiltonian structure on the Lie groups, Lie bialgebras and
the geometric sense of the classical Yang-Baxter equations, Soviet
Math. Dokl. 27 (1983) 68-71.
\item[[8]] M.A. Semonov-Tian-Shansky, What is a classical R-matrix?
Funct. Anal. Appl. 17 (1983) 259-272.
\item[[9]] B.A. Kupershmidt, What a classical $r$-matrix really is, J.
Nonlinear Math. Phy. 6 (1999) 448-488.
\item[[10]] L. Auslander, Simply transitive groups of affine motions.
Amer. J. Math. 99 (1977) 809-826.
\item[[11]] J.-L. Koszul, Domaines born\'es homog\`enes et orbites
de groupes de transformations affines, Bull. Soc. Math. France 89
(1961) 515-533.
\item[[12]] M. Gerstenhaber, The cohomology structure of an
associative ring, Ann. Math. 78 (1963) 267-288.
\item[[13]] E.B. Vinberg, Convex homogeneous cones, Transl. of Moscow Math.
Soc. No. 12 (1963) 340-403.
\item[[14]] A. Andrada, S. Salamon, Complex product structure on Lie
algebras,  Forum Math. 17 (2005) 261-295.
\item[[15]] B.Y. Chu, Symplectic homogeneous spaces, Trans. Amer. Math. Soc.
197 (1974) 145-159.
\item[[16]] J.M. Dardie, A. Medina, Double extension symplectique d'un
groupe de Lie symplectique, Adv. Math. 117 (1996) 208-227.
\item[[17]] J.M. Dardie, A. Medina, Alg\`ebres de Lie k\"ahl\'eriennes et double extension, J. Algebra 185 (1995) 774-795.
\item[[18]] A. Lichnerowicz, A. Medina, On Lie groups with left invariant
symplectic or kahlerian structures, Lett. Math. Phys. 16 (1988)
225-235.
\item[[19]] S.I. Svinolupov, V.V. Sokolov, Vector-matrix generalizations of
classical integrable equations, Theoret. and Math. Phys. 100 (1994)
959-962.
\item[[20]] A. Winterhalder, Linear Nijenhuis-tensors and the construction of
integrable systems, arXiv: physics/9709008.
\item[[21]] A.A. Balinskii, S.P. Novikov, Poisson brackets of hydrodynamic
type, Frobenius algebras and Lie algebras, Soviet Math. Dokl. 32
(1985) 228-231.
\item[[22]] I.M. Gel'fand,  I. Ya. Dorfman, Hamiltonian operators and
algebraic structures related to them, Funct. Anal. Appl. 13 (1979)
248-262.
\item[[23]] E.I. Zel'manov,
On a class of local translation invariant Lie algebras, Soviet Math.
Dokl. 35 (1987) 216-218.
\item[[24]] B. Bakalov, V. Kac, Field algebras, Int. Math. Res. Not.
(2003) 123-159.
\item[[25]] A. Connes, D. Kreimer, Hopf
 algebras, renormalization and noncommutative geometry,
 Comm. Math. Phys. 199 (1998) 203-242.
\item[[26]] F. Chapoton, M. Livernet, Pre-Lie algebras and the rooted trees
operad, Int. Math. Res. Not. (2001) 395-408.
\item[[27]] D. Burde, Left-symmetric algebras, or pre-Lie algebras in
geometry and physics, Cent. Eur. J. Math. 4 (2006) 323-357.
\item[[28]] M. Bordemann, Generalized Lax pairs,
the modified classical Yang-Baxter equation, and affine geometry of
Lie groups, Comm. Math. Phys. 135 (1990) 201-216.
\item[[29]] A. Diatta, A. Medina, Classical Yang-Baxter equation and
left-invariant affine geometry on Lie groups, arXiv:math.DG/0203198.
\item[[30]] P. Etingof, A. Soloviev, Quantization of geometric classical
$r$-matrix, Math. Res. Lett. 6 (1999) 223-228.
\item[[31]] I.Z. Golubschik, V.V. Sokolov, Generalized operator Yang-Baxter
equations, integrable ODEs and nonassociative algebras, J. Nonlinear
Math. Phys., 7 (2000) 184-197.
\item[[32]] N. Jacobson, Lie algebras, Interscience, New York (1962).
\item[[33]] H. Kim, Complete left-invariant affine structures on nilpotent
Lie groups, J. Diff. Geom. 24 (1986) 373-394.
\item[[34]] A. Medina, Flat left-invariant connections adapted to the
automorphism structure of a Lie group, J. Diff. Geom. 16 (1981)
445-474.
\item[[35]] P. Libermann, Sur le probl\`eme d'\'equivalence de certaines
structures infinit\'esimals, Ann. Mat. Pura Appl. 36 (1954) 27-120.
\item[[36]] S. Kaneyuki, Homogeneous symplectic manifolds and
dipolarizations in Lie algebras, Tokyo J. Math. 15 (1992) 313-325.
\item[[37]] B.A. Kupershmidt, Non-abelian phase spaces, J. Phys. A: Math.
Gen. 27 (1994) 2801-2810.
\item[[38]] B.A. Kupershmidt, On the nature of the Virasoro algebra, J.
Nonlinear Math. Phy. 6 (1999) 222-245.
\item[[39]] C.M. Bai, A further study on non-abelian phase
spaces: Left-symmetric algebraic approach and related geometry,
 Rev. Math. Phys. 18 (2006), 545-564.
\item[[40]] C.M. Bai, D.J. Meng, A Lie algebraic approach to Novikov
algebras, J. Geom. Phys. 45 (2003) 218-230.
\item[[41]] N. Jacobson, A note on
automorphisms and derivations of Lie algebras, Proc. Amer. Math.
Soc. 6 (1955) 281-283.
\end{description}

\end{document}